\documentclass[12pt,reqno]{amsart}
\usepackage{amsmath,amsthm,amssymb,amsfonts,amscd}
\usepackage{mathrsfs}
\usepackage{bbding}
\usepackage{graphicx,latexsym}
\usepackage{hyperref}

\newcommand{\bea}{\begin{eqnarray}}
\newcommand{\eea}{\end{eqnarray}}
\newcommand{\bna}{\begin{eqnarray*}}
\newcommand{\ena}{\end{eqnarray*}}

\numberwithin{equation}{section} 

\theoremstyle{plain}

\theoremstyle{definition}

\newtheorem{remark}{Remark}

\begin{document}

\title{Shifted convolution sums of $GL_3$ cusp forms with $\theta$-series}

\author{Qingfeng Sun}
\address{School of Mathematics and Statistics \\ Shandong University, Weihai \\ Weihai \\Shandong 264209 \\China}
\email{qfsun@sdu.edu.cn}
\date{\today}

\begin{abstract}
  Let $A_f(1,n)$ be the normalized Fourier coefficients of a Hecke-Maass cusp form $f$
for $SL_3(\mathbb{Z})$ and
$$
r_3(n)=\#\left\{(n_1,n_2,n_3)\in \mathbb{Z}^3:n_1^2+n_2^2+n_3^2=n\right\}.
$$
Let $1\leq h\leq X$ and $\phi(x)$ be a smooth function compactly
supported on $[1/2,1]$.
It is proved that for any $\varepsilon>0$,
$$
\sum_{n\geq 1}A_f(1,n+h)r_3(n)\phi\left(\frac{n}{X}\right)
\ll_{f,\varepsilon} X^{\frac{3}{2}-\frac{1}{8}+\varepsilon}
$$
uniformly with respect to the shift $h$.
\end{abstract}

\maketitle
\tableofcontents

\section{Introduction}

While the shifted convolution problems involving $GL_2$ automorphic forms are investigated intensively and deep results and
applications have been established (see for example \cite{BH}, \cite{DFI}, \cite{HM}, \cite{H}, \cite{J},
\cite{KMV}, \cite{LLY}, \cite{LS}, \cite{M}, \cite{S}), the situation becomes
much harder for $GL_3$ automorphic forms. For the cube of the Riemann zeta-function $\zeta^3(s)$,
whose coefficients are $\tau_3(n)$, where $\tau_{\ell}(n)=\sum_{n_1n_2\cdot\cdot\cdot n_{\ell}=n}1$,
Pitt \cite{Pitt}, \cite{Pitt2} first proved that, for $r>0$ an integer,
\bea
\sum_{n\leq X}\tau_3(n)a_g(rn-1)\ll_{g,\varepsilon}r^{\frac{5}{24}}X^{\frac{71}{72}+\varepsilon},
\eea
where $a_g(n)$ are the normalized Fourier coefficients of a holomorphic
cusp form for $SL_2(\mathbb{Z})$. Recently, Munshi \cite{Munshi2} improved (1.1) by showing that
$r^{\frac{5}{24}}X^{\frac{71}{72}+\varepsilon}$ can be replaced by
$r^{\frac{2}{7}}X^{\frac{34}{35}+\varepsilon}$. His method is Jutila's version of the
circle method combined with the idea of factorable moduli.
Moreover, let $A_f(1,n)$  and $a_g(n)$ be the normalized Fourier coefficients of a
Hecke-Maass cusp form for $SL_3(\mathbb{Z})$ and $SL_2(\mathbb{Z})$, respectively.
Let $0\leq h\leq X$ and $\varphi(x)$ be a smooth function compactly supported in
$[1,2]$. He also proved by
the same method that (see \cite{M})
\bea
\sum_{n\geq 1}A_f(1,n)a_g(n+h)\varphi\left(\frac{n}{X}\right)\ll_{f,g,\varepsilon}X^{1-\frac{1}{20}+\varepsilon}
\eea
holds uniformly respect to the shift $h$.

Let
$$
r_\ell(n)=\#\left\{(n_1,n_2,\ldots,n_\ell)\in \mathbb{Z}^\ell:n_1^2+\cdots+n_\ell^2=n\right\}.
$$
Then $r_{\ell}(n)$ is the $n$-th Fourier coefficient of the modular form $\theta^\ell(z)$, where
$\theta(z)$ is the classical Jacobi theta series
$$
\theta(z)=\sum_{n\in \mathbb{Z}}e(n^2z).
$$
The shifted convolution sum
\bna
\mathcal {S}_h(X)=\sum_{n\leq X}a_g(n+h)r_{\ell}(n),
\ena
with $a_g(n)$ the normalized Fourier coefficients of
a holomorphic cusp form of weight $\kappa$ for $\Gamma_0(N)$, was first studied by Luo \cite{L}. Precisely,
Luo first established a Voronoi formula for $r_\ell(n)$ and then applied
Poincar\'{e} series reduction and his Voronoi formula to prove that
for $\ell\geq 2$ and $\kappa\geq \frac{\ell}{2}+3$,
\begin{equation}
\mathcal {S}_h(X)\ll_{h,g,\ell,\varepsilon}X^{\frac{\ell}{2}-\vartheta_\ell+\varepsilon},
\end{equation}
where
$$
\vartheta_\ell=\frac{\ell-1}{4(g+1)}, \qquad
g=\left\{\begin{array}{ll}\frac{1+\ell}{2},&\mbox{if $\ell$ odd,}\\
1+\frac{\ell}{2},&\mbox{if $\ell$ even.}\end{array}\right.
$$
In particular,
$
\vartheta_2=\frac{1}{12}, \, \vartheta_3=\frac{1}{6}.
$
Recently, L\"{u}, Wu and Zhai \cite{LWZ} improved Luo's result by the circle method and
showed that (1.3) holds uniformly
for $0\leq h\leq X$ and all $\kappa$, and $\vartheta_\ell$
can be taken as
$
\theta_3=\frac{1}{4}, \,\vartheta_\ell= \frac{1}{2}\, (\ell\geq 4).
$
Moreover, for $N=1$, they proved that
$
\theta_2=\frac{1}{6}, \,\vartheta_\ell= \frac{2}{3}\, (\ell\geq 6).
$

Usually the smaller $\ell$'s are more interesting but the related shifted
convolution problems are more difficult. In this paper, we want to prove the following result.

\medskip

\noindent{\bf Theorem 1.1.}{\it\,
  Let $A_f(1,n)$ be the normalized Fourier
coefficients of a Hecke-Maass cusp form for $SL_3(\mathbb{Z})$ and $1\leq h\leq X$.
Suppose that $\phi(x)$ is a smooth function supported on $[1/2,1]$
satisfying $\phi^{(j)}(x)\ll_j 1$. Then for $X>1$ and any $\varepsilon>0$, we have
\bna
\sum_{n\geq 1}A_f(1,n+h)r_3(n)\phi\left(\frac{n}{X}\right)
\ll_{f,\varepsilon}X^{\frac{3}{2}-\frac{1}{8}+\varepsilon}
\ena
uniformly respect to the shift $h$.}

\begin{remark}
  Under the Ramanujan conjecture, we have $A_f(1,n)\ll_{\varepsilon} n^{\varepsilon}$.
Then by the asymptotic formula in the three-dimensional sphere problem (see \cite{HB2})
\bna
\sum_{n\leq X}r_3(n)=\frac{4\pi}{3}X^{\frac{3}{2}}
+O\left(X^{\frac{21}{32}}\right),
\ena
we have the (conditional) bound
\bea
\sum_{n\geq 1}A_f(1,n+h)r_3(n)\phi\left(\frac{n}{X}\right)\ll_{\varepsilon} X^{\frac{3}{2}+\varepsilon}.
\eea
Further, by the argument in \cite{LWZ} (see Section 3), we have
\bna
\sum_{n\geq 1}A_f(1,n+h)r_3(n)\phi\left(\frac{n}{X}\right)
\ll_{\varepsilon} X^{\frac{3}{2}-\left(\frac{3}{4}-\vartheta\right)+\varepsilon},
\ena
where $\vartheta$ is the exponent in the estimate
\bna
\sum_{n\geq 1}A_f(1,n)e(n \alpha)\phi\left(\frac{n}{X}\right)\ll_{f,\varepsilon} X^{\vartheta+\varepsilon}
\ena
which is uniform in $\alpha \in \mathbb{R}$. Using the best result due to Miller \cite{Miller},
we can take $\vartheta=\frac{3}{4}$ and also obtain the same bound as in (1.4).
The bound in Theorem 1.1 is unconditional and better than (1.4) by
a power of $X^{\frac{1}{8}-\varepsilon}$.
\end{remark}

\begin{remark}
 Although the shifted convolution sum in question should be comparable
to the sum in (1.2) in view of the Voronoi formula for $r_{\ell}(n)$ established recently by Luo \cite{L},
the methods in Munshi \cite{M} which obtains nontrivial power saving from the structural advantage
of Jutila's variation of the circle method, are not applicable
in our situation. This is because that there is
no good upper bound for the exponential twisted sum
\bna
\sum_{n\leq X}r_{\ell}(n)e(n \alpha), \qquad \alpha \in \mathbb{R}
\ena
which is, however, necessary for Jutila's version of the circle method.
\end{remark}

\begin{remark} For $h=0$, the bound in Theorem 1.1 still holds.
Moreover, Theorem 1.1 can be generalized to the cases
$r_{\ell}(n)$, $\ell\geq 4$ without any difficulty.
It is also worth noting that replacing $A_f(1,n)$ by the
triple divisor function $\tau_3(n)$,
we can get a better exponent saving
than $\frac{1}{8}-\varepsilon$ due to the bound $\tau_3(n)\ll n^{\varepsilon}$
for any $\varepsilon>0$. In fact, we have
\bna
\sum_{n\geq 1}\tau_3(n+h)r_3(n)\phi\left(\frac{n}{X}\right)
=4\mathcal {C}_2\mathcal{I}_0X^{\frac{3}{2}}+
4\mathcal {C}_1\mathcal{I}_1X^{\frac{3}{2}}+2\mathcal {C}_0\mathcal{I}_2X^{\frac{3}{2}}
+O_{\varepsilon}\left(X^{\frac{3}{2}-\frac{1}{4}+\varepsilon}\right),
\ena
where for $\ell=0,1,2$,
\bna
\mathcal {C}_{\ell}:=\mathcal {C}_{\ell}(h)=\sum_{q=1}^{\infty}\frac{1}{q^5}\sum_{n|q}n\tau(n)P_{\ell}(n,q)
\sum_{a=1\atop (a,q)=1}^qe\left(\frac{ah}{q}\right)
G(a,0;q)^3 S\left(-\overline{a},0;\frac{q}{n}\right),
\ena
and
\bna
\mathcal{I}_{\ell}:=\mathcal{I}_{\ell}(X,h,\phi)=\int_{-\infty}^{\infty}
\left(\int_0^{\infty}\phi\left(\frac{u}{X}\right)e(-\beta u)
(\log(u+h))^{\ell}\mathrm{d}u\right)
\left(\int_0^{1}e(\beta Xv^2)\mathrm{d}v\right)^3\mathrm{d}\beta.
\ena
Here $\tau(n)=\sum_{d|n}1$ is the divisor function, $\overline{a}$ denotes the multiplicative inverse of $a\bmod q$, $S(a,b;c)$ is the classical Kloosterman sum, $G(a,b;q)$ is the Gauss sum
\bna
G(a,b;q)=\sum\limits_{d\bmod q}e\left(\frac{ad^2+bd}{q}\right),
\ena
$P_0(n,q)=1$ and $P_\ell(n,q)$ ($\ell=1,2$) are given by
\bna
P_1(n,q)&=&\frac{5}{3}\log n-3\log q+3\gamma-\frac{1}{3\tau(n)}
\sum_{d|n}\log d,\\
P_2(n,q)&=&\left(\log n\right)^2-5\log q \log n
+\frac{9}{2}(\log q)^2+3\gamma^2-3\gamma_1+7\gamma \log n-9\gamma \log q\nonumber\\
&&+\frac{1}{\tau(n)}\left(\left(\log n+\log q-5\gamma\right)\sum_{d|n}\log d-\frac{3}{2}
\sum_{d|n}(\log d)^2\right)
\ena
with $\gamma:=\lim\limits_{s\rightarrow 1}\left(\zeta(s)-\frac{1}{s-1}\right)$
being the Euler constant and $\gamma_1:=-\frac{\mathrm{d}}{\mathrm{d}s}\left.\left(\zeta(s)-\frac{1}{s-1}\right)\right|_{s=1}$
being the Stieltjes constant. The proof is very similar as that
of Theorem 1.1 (see also \cite{SZ}).

\end{remark}

Let us briefly outline the method of the proof.
Applying the classical Hardy-Littlewood-Kloosterman circle method to
the shifted convolution sum in Theorem 1.1, we arrive at an expression which is convenient
to apply the Voronoi formula for $GL_3$ and an asymptotic formula for the sum
$\sum_{|n|\leq X}e(\alpha n^2)$, $\alpha\in [0,1]$.
After that we are left with eight sums (in principle, of same difficulty) involving
higher-dimensional character sums.
In fact, we need to consider the twisted exponential sum (see Section 5)
\bna
\widetilde{\mathscr{T}}(p)=\sum_{x \in \mathbb{F}_p^{\times}}
\left(\frac{x}{p}\right)
e\left(\frac{r_1hx-\overline{4}(4v+b_1^2+b_2^2+b_3^2)\overline{x}}{p}\right)
S(-r_2\overline{x},r_3n;p),
\ena
where $p$ is an odd prime, $(r_i,p)=1$, $i=1,2,3$, and $b_1,b_2,b_3,v,n\in \mathbb{Z}$.
This type of character sums have been studied in the
work of Adolphson and Sperber \cite{AS}.
The main saving in Theorem 1.1 comes from the nontrivial estimation of these
character sums.

\medskip
\noindent
{\bf Notation.}
Throughout the paper, the letters $h$, $q$, $m$ and $n$, with or without subscript,
denote integers. The letter $\varepsilon$ is an arbitrarily small
positive constant, not necessarily the same at different occurrences. The symbol
$\ll_{a,b,c}$ denotes that the implied constant depends at most on $a$, $b$ and $c$.

\section{Preliminaries on $GL_3$ cusp forms}
\setcounter{equation}{0}
\bigskip

Let $f$ be a Hecke-Maass cusp form of type $\nu=(\nu_1,\nu_2)$ for $SL_3(\mathbb{Z})$, normalized so that the first Fourier
coefficient is 1. Then $f$ has a Fourier-Whittaker expansion (see \cite{G})
\bna
f(z)=\sum_{\gamma\in U_2(\mathbb{Z})\backslash SL_2(\mathbb{Z})}
\sum_{n_1=1}^{\infty}\sum_{n_2\neq 0}\frac{A_f(n_1,n_2)}{n_1|n_2|}
W_J\left(M\left(\begin{array}{ll}\gamma &\, \\
\,& 1\end{array}\right)z,\nu,\psi_{1,1}\right)
\ena
where $U_2(\mathbb{Z})$ is the group of $2\times 2$ upper triangular matrices with
integer entries and ones on the diagonal, $W_J\left(z,\nu,\psi_{1,1}\right)$ is the Jacquet-Whittaker function and
$M=\mathrm{diag}(n_1|n_2|,n_1,1)$. By Rankin-Selberg theory,
the Fourier coefficients $A_f(n_1,n_2)$ satisfy
\bea
\sum_{n_2\leq N} \left|A_f(n_1,n_2)\right|\ll_f N|n_1|, \qquad
\sum_{n_1\leq N} \left|A_f(n_1,n_2)\right|\ll_f N|n_2|.
\eea

Let
\bna
\mu_1=-\nu_1-2\nu_2+1, \qquad \mu_2=-\nu_1+\nu_2,\qquad
\mu_3=2\nu_1+\nu_2-1.
\ena
The generalized Ramanujan conjecture asserts that $\mathrm{Re}(\mu_j)=0$, $1\leq j\leq 3$,
while the current record bound due to
Luo, Rudnick and Sarnak \cite{LRS} is
\bea
|\mathrm{Re}(\mu_j)| \leq \frac{1}{2}-\frac{1}{10}, \quad 1\leq j\leq 3.
\eea
Let $\phi(x)$ be a smooth function compactly supported on $(0,\infty)$ and
denote by $\widetilde{\phi}(s)$
the Mellin transform of $\phi(x)$.
For $k=0,1$, we define
\bea
\Phi_k(x):=\int\limits_{\mathrm{Re}(s)=\sigma}(\pi^3 x)^{-s}\prod_{j=1}^3
\frac{\Gamma\left(\frac{1+s+\mu_j+2k}{2}\right)}
{\Gamma\left(\frac{-s-\mu_j}{2}\right)}
\widetilde{\phi}(-s-k)\mathrm{d}s
\eea
with $\sigma>\max\limits_{1\leq j\leq 3}\{-1-\mathrm{Re}(\mu_j)-2k\}$. Set
\bea
\Phi^{\pm}(x)=\Phi_0(x)\pm\frac{1}{i\pi^3x}\Phi_1(x).
\eea
Then we have the following Voronoi-type formula (see \cite{GL1}, \cite{MS}).

\medskip

\noindent {\bf Lemma 2.1.} {\it Let $A_f(m,n)$ be the
Fourier coefficients of a Maass cusp form for $SL_3(\mathbb{Z})$.
Suppose that $\phi(x) \in C_c^\infty(0,\infty)$. Let
$a,q\in \mathbb{Z}$ with $q\geq 1$, $(a,q)=1$ and $a \overline{a} \equiv 1(\bmod q)$. Then
\bna
\sum_{n\geq 1}A_f(m,n)e\left(\frac{an}{q}\right)\phi(n)&=&\frac{q\pi^{-\frac{5}{2}}}{4i}\sum_{n_1|qm}
\sum_{n_2=1}^{\infty}
\frac{A_f(n_2,n_1)}{n_1n_2}S(m\overline{a},n_2;mqn_1^{-1})\Phi^+\left(\frac{n_1^2n_2}{q^3m}\right)\\ &+&
\frac{q\pi^{-\frac{5}{2}}}{4i}\sum_{n_1|qm}\sum_{n_2=1}^{\infty}
\frac{A_f(n_1,n_2)}{n_1n_2}S(m\overline{a},-n_2;mqn_1^{-1})\Phi^-\left(\frac{n_1^2n_2}{q^3m}\right),
\ena
where $S(m,n;c)$ is the classical Kloosterman sum.
}

The functions $\Phi^{\pm}(x)$ has the following properties.

\medskip

\noindent {\bf Lemma 2.2.} {\it Suppose that $\phi(x)$ is a smooth function of
compact support in $[AX,BX]$,
where $X>0$ and $B>A>0$,
satisfying $\phi^{(j)}(x)\ll_{A,B,j}P^j $ for any integer $j\geq 0$. Then for $x>0$ and any integer $N\geq 0$, we have
$$
\Phi^{\pm}(x)\ll_{f,A,B,N,\varepsilon}(xX)^{-\varepsilon}(PX)^3\left(\frac{x}{P^3X^2}\right)^{-N}.
$$
}

\noindent{\bf Proof.} Let $\Phi_k(x)$ be as in (2.3). Changing variable $s+k\rightarrow s$, we have,
for $\sigma>\max\limits_{1\leq j\leq 3}\{-1-\mathrm{Re}(\mu_j)-k\}$,
\bea
\Phi_k(x)=(\pi^3 x)^k\int\limits_{\mathrm{Re}(s)=\sigma}(\pi^3 x)^{-s}\prod_{j=1}^3
\frac{\Gamma\left(\frac{1+s+\mu_j+k}{2}\right)}
{\Gamma\left(\frac{-s-\mu_j+k}{2}\right)}
\widetilde{\phi}(-s)\mathrm{d}s,
\eea
where
\bna
\widetilde{\phi}(s)=\int_0^{\infty}\phi(u)u^{s-1}\mathrm{d}u.
\ena
Let $s=\sigma+i t$. By partial integration $j$ times, we have, for $\sigma\not\in \mathbb{Z}$,
\bea
\widetilde{\phi}(-s)=\frac{1}{s(s-1)\cdot\cdot\cdot (s-j+1)}
\int_0^{\infty}\phi^{(j)}(u)u^{-s+j-1}\mathrm{d}u \ll_{A,B,\sigma,j} X^{-\sigma}(PX)^j(1+|t|)^{-j}.
\eea
Moreover, by Stirling's formula,
\bea
\prod_{j=1}^3
\frac{\Gamma\left(\frac{1+s+\mu_j+k}{2}\right)}
{\Gamma\left(\frac{-s-\mu_j+k}{2}\right)}\ll_{f,\sigma}\prod_{j=1}^3(1+|t|)^{\sigma+\frac{1}{2}+\mathrm{Re}(\mu_j)}
\ll_{f,\sigma}(1+|t|)^{3\sigma+\frac{3}{2}}
\eea
in view of the fact that $\sum_{j=1}^3\mu_j=0$.

Moving the contour to $\sigma=N+\varepsilon$ with $N\in \mathbb{Z}^+\cup 0$ and taking $j=3N+3$,
by (2.5)-(2.7), we have
\bea
\Phi_k(x)&\ll_{f,A,B,\sigma,j}&x^k(xX)^{-\sigma}(PX)^j\int_{-\infty}^{+\infty}
(1+|t|)^{-j+3\sigma+\frac{3}{2}}\mathrm{d}t\nonumber\\
&\ll_{f,A,B,N,\varepsilon}&x^k(xX)^{-N-\varepsilon}(PX)^{3N+3}
\int_{-\infty}^{+\infty}
(1+|t|)^{-\frac{3}{2}+3\varepsilon}\mathrm{d}t\nonumber\\
&\ll_{f,A,B,N,\varepsilon}&x^k(xX)^{-\varepsilon}(PX)^3\left(\frac{x}{P^3X^2}\right)^{-N}.
\eea
Then Lemma 2.2 follows from (2.4) and (2.8). \hfill $\Box$

\medskip

By Lemma 2.2, for any fixed $\varepsilon>0$ and $xX\geq X^{\varepsilon}(PX)^3$,
$\Phi^{\pm}(x)$ is negligibly small. For $xX\ll X^{\varepsilon}$, as in \cite{Li},
we move the line of integration in (2.6) to $\sigma=-\frac{11}{20}$ (guaranteed by (2.2)),
by (2.6) with $j=1$ and (2.7), we have
\bna
\Phi_k(x)\ll_{f,A,B}x^k(xX)^{\frac{11}{20}}(PX)
\int_{-\infty}^{+\infty}
(1+|t|)^{-\frac{23}{20}}\mathrm{d}t
\ll_{f,A,B,\varepsilon}x^kPX^{1+\varepsilon},
\ena
and by (2.4),
\bea
\Phi^{\pm}(x)\ll_{f,A,B,\varepsilon}PX^{1+\varepsilon}.
\eea

\medskip

For $xX\gg X^{\varepsilon}$, we shall use the following result (see \cite{Li1}, \cite{RY}).
For $\mu_1=\mu_2=\mu_3=0$, this was proved in \cite{Iv}.

\medskip

\noindent {\bf Lemma 2.3.} {\it
Suppose that $\phi(x)$ is a smooth function of compact support on $[AX,BX]$,
where $X>0$ and $B>A>0$. Then for $x>0$, $xX\gg 1$, $\ell\geq 2$ and $k=0,1$, we have
\bna
\Phi_k(x)&=&(\pi^3 x)^{k+1}\sum_{j=1}^{\ell} \int_0^{\infty}\phi(u)
\left(a_k(j)e\left(3(x u)^{\frac{1}{3}}\right)+b_k(j)e\left(-3(x u)^{\frac{1}{3}}\right)\right)
\frac{\mathrm{d}u}{(\pi^3xu)^{\frac{j}{3}}}\\
&&+O_{f,A,B,\varepsilon,\ell}\left((\pi^3x)^k(\pi^3xX)^{-\frac{\ell}{3}+\frac{1}{2}+\varepsilon}\right),
\ena
where $a_k(j)$, $b_k(j)$ are constants with
\bna
a_0(1)=-\frac{2\sqrt{3\pi}}{3}, \quad b_0(1)=\frac{2\sqrt{3\pi}}{3},
\quad a_0(1)=b_1(1)=-\frac{2\sqrt{3\pi}}{3}i.
\ena
}

\medskip

\section{Proof of Theorem 1.1}
\setcounter{equation}{0}
\medskip

Denote
\bna
\mathscr{S}_h(X)=\sum_{n\geq 1}A_f(1,n+h)r_3(n)\phi\left(\frac{n}{X}\right).
\ena
We first transform $\mathscr{S}_h(X)$ by the Hardy-Littlewood-Kloosterman circle method
(see for example, \cite{Iwan}, Section 11.4).
Let
\bna
\mathscr{F}(\alpha)=\sum_{|m|\leq \sqrt{X}}e(\alpha m^2)
\ena
and
\bea
\mathscr{G}(\alpha)=\sum_{n\geq 1}A_f(1,n+h)e(-\alpha n)\phi\left(\frac{n}{X}\right).
\eea
Then $\mathscr{S}_h(X)$ can be written as
\bna
\mathscr{S}_h(X)=\int_0^1\mathscr{F}^3(\alpha)\mathscr{G}(\alpha)\mathrm{d}\alpha.
\ena
Note that $\mathscr{F}^3(\alpha)\mathscr{G}(\alpha)$ is a periodic function of period 1. We have
\bna
\mathscr{S}_h(X)=\int_{-1/(Q+1)}^{1-1/(Q+1)}\mathscr{F}^3(\alpha)\mathscr{G}(\alpha)\mathrm{d}\alpha,
\ena
where $Q=[5\sqrt{X}]$.
Dissecting the unit interval with Farey's points of order $Q$, we have
\bna
\mathscr{S}_h(X)=\sum_{q\leq Q}\sideset{}{^*}\sum_{a=1}^q
\int\limits_{\mathscr{M}(a,q)}\mathscr{F}^3\left(\frac{a}{q}+\beta\right)\mathscr{G}\left(\frac{a}{q}+\beta\right)\mathrm{d}\beta,
\ena
where the $*$ denotes the condition $(a,q)=1$,
\[
\mathscr{M}(a,q)=\left[-\frac{1}{q(q+q')},\frac{1}{q(q+q'')}\right],
\]
$\frac{a'}{q'}$, $\frac{a}{q}$ and $\frac{a''}{q''}$ are consecutive Farey fractions and they are determined
by the conditions
\[
Q<q+q',q+q''\leq q+Q, \quad aq'\equiv 1(\bmod q), \quad aq''\equiv -1(\bmod q).
\]
Exchanging the order of the summation over $a$ and the integration over $\beta$ as in
Heath-Brown \cite{HB}, we have
\bea
\mathscr{S}_h(X)=\sum_{q\leq Q}\int\limits_{|\beta|\leq \frac{1}{qQ}}\sum_{v\bmod q}
\varrho(v,q,\beta)\sideset{}{^*}\sum_{a=1}^qe\left(-\frac{\overline{a}v}{q}\right)
\mathscr{F}^3\left(\frac{a}{q}+\beta\right)\mathscr{G}\left(\frac{a}{q}+\beta\right)\mathrm{d}\beta,
\eea
where $\varrho(v,q,\beta)$ satisfies
\bea
\varrho(v,q,\beta)\ll \frac{1}{1+|v|}.
\eea

For an asymptotic formula of $\mathscr{F}\left(\frac{a}{q}+\beta\right)$,
we quote the following result (see Theorem 4.1 in \cite{V} or Lemma 4.1 in \cite{Z}).

\medskip

\noindent {\bf Lemma 3.1.} {\it Suppose that $(a,q)=1$, $q\leq Q$ and
$|\beta|\leq 1/(qQ)$. We have
\bea
\mathscr{F}\left(\frac{a}{q}+\beta\right)=\frac{2G(a,0;q)}{q}\Psi_0(\beta)
+\sum_{-\frac{3q}{2}<b\leq \frac{3q}{2}}G(a,b;q)\Psi(b,q,\beta),
\eea
where $G(a,b;q)$ is the Gauss sum
\bea
G(a,b;q)=\sum\limits_{x\bmod q}e\left(\frac{ax^2+bx}{q}\right),
\eea
$\Psi_0(\beta)$ is the integral
\bea
\Psi_0(\beta)=\int_0^{\sqrt{X}}e(\beta x^2)\mathrm{d}x,
\eea
and $\Psi(b,q,\beta)$ satisfies
\bea
\sum_{-\frac{3q}{2}<b\leq \frac{3q}{2}}|\Psi(b,q,\beta)|\ll \log(q+2).
\eea
}

For $\mathscr{G}(\alpha)$ in (3.1), we apply Lemma 2.1
with $\phi_{\beta}(x)=\phi\left(\frac{x-h}{X}\right)e(-\beta x)$ getting
\bea
\mathscr{G}\left(\frac{a}{q}+\beta\right)
&=&e\left(\frac{ha}{q}+h \beta\right)\sum_{n\geq 1}A_f(1,n)e\left(-\frac{a n}{q}\right)
\phi\left(\frac{n-h}{X}\right)e(-\beta n)\nonumber\\
&=&\frac{q\pi^{-\frac{5}{2}}}{4i}e\left(\frac{ha}{q}+h \beta\right)
\sum_{n_1|q}\sum_{n_2=1}^{\infty}\frac{A_f(n_2,n_1)}{n_1n_2}S\left(-\overline{a},n_2;\frac{q}{n_1}\right)
\Phi_{\beta}^+\left(\frac{n_1^2n_2}{q^3}\right)\nonumber\\
&&+\frac{q\pi^{-\frac{5}{2}}}{4i}e\left(\frac{ha}{q}+h \beta\right)
\sum_{n_1|q}\sum_{n_2=1}^{\infty}\frac{A_f(n_1,n_2)}{n_1n_2}S\left(-\overline{a},-n_2;\frac{q}{n_1}\right)
\Phi_{\beta}^-\left(\frac{n_1^2n_2}{q^3}\right),\nonumber\\
\eea
where
\bea
\Phi_{\beta}^{\pm}(x)=\Phi_0(x,\beta)\pm\frac{1}{i\pi^3x}\Phi_1(x,\beta)
\eea
with
\bea
\Phi_k(x,\beta)=\int\limits_{\mathrm{Re}(s)=\sigma}(\pi^3 x)^{-s}\prod_{j=1}^3
\frac{\Gamma\left(\frac{1+s+\mu_j+2k}{2}\right)}
{\Gamma\left(\frac{-s-\mu_j}{2}\right)}
\widetilde{\phi_{\beta}}(-s-k)\mathrm{d}s.
\eea

By (3.4) and (3.8), we have
\bea
\sideset{}{^*}\sum_{a=1}^qe\left(-\frac{\overline{a}v}{q}\right)
\mathscr{F}^3\left(\frac{a}{q}+\beta\right)\mathscr{G}\left(\frac{a}{q}+\beta\right)
=\sum_{j=1}^8\mathscr{D}_j(v,q,\beta),
\eea
where
\bea
\mathscr{D}_1(v,q,\beta)&=&\frac{2}{\pi^{\frac{5}{2}}i}\frac{e(h\beta)\Psi_0^3(\beta)}{q^2}
\sum_{n_1|q}\sum_{n_2=1}^{\infty}\frac{A_f(n_2,n_1)}{n_1n_2}\Phi_{\beta}^+\left(\frac{n_1^2n_2}{q^3}\right)
\mathscr{C}(0,0,0,n_1,n_2,h,v;q),\\
\mathscr{D}_2(v,q,\beta)&=&\frac{2}{\pi^{\frac{5}{2}}i}\frac{e(h\beta)\Psi_0^3(\beta)}{q^2}
\sum_{n_1|q}\sum_{n_2=1}^{\infty}\frac{A_f(n_1,n_2)}{n_1n_2}\Phi_{\beta}^-\left(\frac{n_1^2n_2}{q^3}\right)\nonumber\\
&&\times\mathscr{C}(0,0,0,n_1,-n_2,h,v;q),\nonumber\\
\mathscr{D}_3(v,q,\beta)&=&\frac{3}{\pi^{\frac{5}{2}}i}\frac{e(h\beta)\Psi_0^2(\beta)}{q}
\sum_{n_1|q}\sum_{n_2=1}^{\infty}\frac{A_f(n_2,n_1)}{n_1n_2}\Phi_{\beta}^+\left(\frac{n_1^2n_2}{q^3}\right)\nonumber\\
&&\times\sum_{-\frac{3q}{2}<b\leq \frac{3q}{2}}\Psi(b,q,\beta)
\mathscr{C}(0,0,b,n_1,n_2,h,v;q),\\
\mathscr{D}_4(v,q,\beta)&=&\frac{3}{\pi^{\frac{5}{2}}i}\frac{e(h\beta)\Psi_0^2(\beta)}{q}
\sum_{n_1|q}\sum_{n_2=1}^{\infty}\frac{A_f(n_1,n_2)}{n_1n_2}\Phi_{\beta}^-\left(\frac{n_1^2n_2}{q^3}\right)\nonumber\\
&&\times\sum_{-\frac{3q}{2}<b\leq \frac{3q}{2}}\Psi(b,q,\beta)
\mathscr{C}(0,0,b,n_1,-n_2,h,v;q),\nonumber\\
\mathscr{D}_5(v,q,\beta)&=&\frac{3}{2\pi^{\frac{5}{2}}i}e(h\beta)\Psi_0(\beta)
\sum_{n_1|q}\sum_{n_2=1}^{\infty}\frac{A_f(n_2,n_1)}{n_1n_2}\Phi_{\beta}^+\left(\frac{n_1^2n_2}{q^3}\right)
\sum_{-\frac{3q}{2}<b_j\leq \frac{3q}{2}\atop j=1,2}\Psi(b_1,q,\beta)\Psi(b_2,q,\beta)\nonumber\\
&&\times\mathscr{C}(0,b_1,b_2,n_1,n_2,h,v;q),\\
\mathscr{D}_6(v,q,\beta)&=&\frac{3}{2\pi^{\frac{5}{2}}i}e(h\beta)\Psi_0(\beta)
\sum_{n_1|q}\sum_{n_2=1}^{\infty}\frac{A_f(n_1,n_2)}{n_1n_2}\Phi_{\beta}^-\left(\frac{n_1^2n_2}{q^3}\right)
\sum_{-\frac{3q}{2}<b_j\leq \frac{3q}{2}\atop j=1,2}\Psi(b_1,q,\beta)\Psi(b_2,q,\beta)\nonumber\\
&&\times\mathscr{C}(0,b_1,b_2,n_1,-n_2,h,v;q),\nonumber
\eea
\bea
\mathscr{D}_7(v,q,\beta)&=&\frac{q}{4\pi^{\frac{5}{2}}i}e(h\beta)
\sum_{n_1|q}\sum_{n_2=1}^{\infty}\frac{A_f(n_2,n_1)}{n_1n_2}\Phi_{\beta}^+\left(\frac{n_1^2n_2}{q^3}\right)
\sum_{-\frac{3q}{2}<b_j\leq \frac{3q}{2}\atop 1\leq j\leq 3}
\Psi(b_1,q,\beta)\Psi(b_2,q,\beta)\Psi(b_3,q,\beta)\nonumber\\
&&\times\mathscr{C}(b_1,b_2,b_3,n_1,n_2,h,v;q),\\
\mathscr{D}_8(v,q,\beta)&=&\frac{q}{4\pi^{\frac{5}{2}}i}e(h\beta)
\sum_{n_1|q}\sum_{n_2=1}^{\infty}\frac{A_f(n_1,n_2)}{n_1n_2}\Phi_{\beta}^-\left(\frac{n_1^2n_2}{q^3}\right)
\sum_{-\frac{3q}{2}<b_j\leq \frac{3q}{2}\atop 1\leq j\leq 3}
\Psi(b_1,q,\beta)\Psi(b_2,q,\beta)\Psi(b_3,q,\beta)\nonumber\\
&&\times\mathscr{C}(b_1,b_2,b_3,n_1,-n_2,h,v;q)\nonumber
\eea
with
\bea
\mathscr{C}(b_1,b_2,b_3,n_1,n_2,h,v;q)=\sideset{}{^*}\sum_{a=1}^qe\left(\frac{ah-\overline{a}v}{q}\right)
G(a,b_1;q)G(a,b_2;q)G(a,b_3;q) S\left(-\overline{a},n_2;\frac{q}{n_1}\right).
\nonumber\\
\eea
We only estimate the contributions from $\mathscr{D}_1$, $\mathscr{D}_3$,
$\mathscr{D}_5$, $\mathscr{D}_7$, and the contributions from $\mathscr{D}_2$, $\mathscr{D}_4$,
$\mathscr{D}_6$, $\mathscr{D}_8$ can be estimated similarly.

The following propositions will be proved in the next section.

\medskip

\noindent {\bf Proposition 3.2.} {\it For any $\varepsilon>0$, we have
\bea
\sum_{n_2=1}^{\infty}\frac{|A_f(n_2,n_1)|}{n_1n_2}
\left|\Phi_{\beta}^{\pm}\left(\frac{n_1^2n_2}{q^3}\right)\right|\ll_{f,\varepsilon}
X^{\varepsilon}(1+|\beta|X)^2
\eea
and
\bea
\sum_{n_2=1}^{\infty}\frac{|A_f(n_1,n_2)|}{n_1n_2}
\left|\Phi_{\beta}^{\pm}\left(\frac{n_1^2n_2}{q^3}\right)\right|
\ll_{f,\varepsilon}
X^{\varepsilon}(1+|\beta|X)^2.
\eea
}

\noindent {\bf Proposition 3.3.} {\it
Let $q_1$ be the largest factor of $q$ such that $q_1|n_1$ and $(q_1, q/q_1)=1$.
Let $q_2$ be the largest factor of $q/q_1$ such that $q_2|n_1^{\infty}$ and
$\left(q_2,\frac{q}{q_1q_2}\right)=1$.
Let $q=q_1q_2q_3'q_3''$, $(q_3',2q_3'')=1$, $q_3'$ square-free and $4q_3''$ square-full.
Then for any $\varepsilon>0$, we have
\bna
\mathscr{C}(b_1,b_2,b_3,n_1,n_2,h,v;q)\ll_{\varepsilon} \frac
{(q_1q_2q_3'')^{3+\varepsilon}q_3'^{\frac{5}{2}+\varepsilon}(h,q_3')^{\frac{1}{2}}}{\sqrt{n_1}}.
\ena
}

By the second derivative test and the trivial estimation,
$\Psi_0(\beta)$ in (3.6) is bounded by
\bea
\Psi_0(\beta)\ll \left(\frac{X}{1+|\beta|X}\right)^{\frac{1}{2}}.
\eea
Note that the condition $\frac{n_1^2n_2}{q^3}X< X^{\varepsilon}(1+|\beta|X)^3$ with
$|\beta|\leq 1/(qQ)$ implies that $n_1^2n_2\ll X^{2+\varepsilon}/Q^3\ll X^{\frac{1}{2}+\varepsilon}$.

Let $q$ be as in Proposition 3.3. Denote $q_0=q_2q_3''$.
Note that $q_0$ is square-full.
By (3.12), (3.17), (3.19), Lemma 2.2 and Proposition 3.3, we have
\bna
\mathscr{D}_1(v,q,\beta)
&\ll_{f,\varepsilon}&
 \frac{1}{q^2}\left(\frac{X}{1+|\beta|X}\right)^{\frac{3}{2}}
\sum_{n_1|q}\sum_{\frac{n_1^2n_2}{q^3}X< X^{\varepsilon}(1+|\beta|X)^3}\frac{\left|A_f(n_2,n_1)\right|}{n_1n_2}
\left|\Phi_{\beta}^+\left(\frac{n_1^2n_2}{q^3}\right)\right|\\
&&\times|\mathscr{C}(0,0,0,n_1,n_2,h,v;q)|+1\nonumber\\
&\ll_{f,\varepsilon}&X^{\varepsilon}
\left(\frac{X}{1+|\beta|X}\right)^{\frac{3}{2}}
\sum_{n_1\ll X^{\frac{1}{4}+\varepsilon}\atop n_1\equiv 0(\bmod q_1)}
\frac{q_1q_0q_3'^{\frac{1}{2}}(h,q_3')^{\frac{1}{2}}}{\sqrt{n_1}}\sum_{n_2=1}^{\infty}\frac{\left|A_f(n_2,n_1)\right|}{n_1n_2}
\left|\Phi_{\beta}^+\left(\frac{n_1^2n_2}{q^3}\right)\right|
+1\nonumber\\
&\ll_{f,\varepsilon}&X^{\frac{3}{2}+\varepsilon}
\left(1+|\beta|X\right)^{\frac{1}{2}}
\sum_{n_1\ll X^{\frac{1}{4}+\varepsilon}\atop n_1\equiv 0(\bmod q_1)}
\frac{q_1q_0q_3'^{\frac{1}{2}}(h,q_3')^{\frac{1}{2}}}{\sqrt{n_1}}.
\ena
By (3.3) and the estimate above, we have
\bea
&&\sum_{q\leq Q}\int\limits_{|\beta|\leq \frac{1}{q Q}}\sum_{v\bmod q}
\varrho(v,q,\beta)\mathscr{D}_1(v,q,\beta)\mathrm{d}\beta\nonumber\\
&\ll_{f,\varepsilon}&X^{\frac{3}{2}+\varepsilon}\sum_{n_1\ll X^{\frac{1}{4}+\varepsilon}}n_1^{-\frac{1}{2}}
\sum_{q_1|n_1}\sum_{q_3'\leq Q/q_1 \atop q_3' \,\mathrm{square-free}}
\sum_{q_0\leq Q/(q_1q_3') \atop 4q_0 \, \mathrm{square-full}}
q_1q_0q_3'^{\frac{1}{2}}(h,q_3')^{\frac{1}{2}}
\int\limits_{|\beta|\leq \frac{1}{q_1q_0q_3'Q}}
(1+|\beta|X)^{\frac{1}{2}}\mathrm{d}\beta\nonumber\\
\nonumber\\
&\ll_{f,\varepsilon}&X^{\frac{3}{2}+\varepsilon}\sum_{n_1\ll X^{\frac{1}{4}+\varepsilon}}n_1^{-\frac{1}{2}}
\sum_{q_1|n_1}\sum_{q_3'\leq Q/q_1 \atop q_3' \,\mathrm{square-free}}
\sum_{q_0\leq Q/(q_1q_3') \atop 4q_0 \, \mathrm{square-full}} q_1q_0q_3'^{\frac{1}{2}}(h,q_3')^{\frac{1}{2}}
 \left(\frac{1}{q_1q_0q_3'Q}+\frac{X^{\frac{1}{2}}}{(q_1q_0q_3'Q)^{\frac{3}{2}}} \right)\nonumber\\
&\ll_{f,\varepsilon}&\frac{X^{\frac{3}{2}+\varepsilon}}{Q}
\sum_{n_1\ll X^{\frac{1}{4}+\varepsilon}}n_1^{-\frac{1}{2}}
\sum_{q_1|n_1}\sum_{q_3'\leq Q/q_1}(h,q_3')^{\frac{1}{2}}q_3'^{-\frac{1}{2}}
\sum_{q_0\leq Q/(q_1q_3')\atop 4q_0\,\mathrm{square-full}}1\nonumber\\
&&+\frac{X^{2+\varepsilon}}{Q^{\frac{3}{2}}}
\sum_{n_1\ll X^{\frac{1}{4}+\varepsilon}}n_1^{-\frac{1}{2}}
\sum_{q_1|n_1}q_1^{-\frac{1}{2}}\sum_{q_3'\leq Q/q_1}(h,q_3')^{\frac{1}{2}}q_3'^{-1}
\sum_{q_0\leq Q/(q_1q_3')\atop 4q_0\,\mathrm{square-full}}q_0^{-\frac{1}{2}}\nonumber\\
&\ll_{f,\varepsilon}&\frac{X^{\frac{3}{2}+\varepsilon}}{Q}
\sum_{n_1\ll X^{\frac{1}{4}+\varepsilon}}n_1^{-\frac{1}{2}}
\sum_{q_1|n_1}\sum_{q_3'\leq Q/q_1}(h,q_3')^{\frac{1}{2}}q_3'^{-\frac{1}{2}}
\left(\frac{Q}{q_1q_3'}\right)^{\frac{1}{2}}\nonumber\\
&&+\frac{X^{2+\varepsilon}}{Q^{\frac{3}{2}}}
\sum_{n_1\ll X^{\frac{1}{4}+\varepsilon}}n_1^{-\frac{1}{2}}
\sum_{q_1|n_1}q_1^{-\frac{1}{2}}\sum_{q_3'\leq Q/q_1}(h,q_3')^{\frac{1}{2}}q_3'^{-1}\nonumber\\
&\ll_{f,\varepsilon}&\frac{X^{\frac{3}{2}+\varepsilon}}{Q^{\frac{1}{2}}}
\sum_{n_1\ll X^{\frac{1}{4}+\varepsilon}}n_1^{-\frac{1}{2}}
\sum_{q_1|n_1}q_1^{-\frac{1}{2}}\sum_{d|h}d^{-\frac{1}{2}}
\sum_{q_4\leq Q/{q_1d}}q_4^{-1}\nonumber\\
&&+\frac{X^{2+\varepsilon}}{Q^{\frac{3}{2}}}
\sum_{n_1\ll X^{\frac{1}{4}+\varepsilon}}n_1^{-\frac{1}{2}}
\sum_{q_1|n_1}q_1^{-\frac{1}{2}}\sum_{d|h}d^{-\frac{1}{2}}\sum_{q_4\leq Q/(q_1d)}q_4^{-1}\nonumber\\
&\ll_{f,\varepsilon}&\frac{X^{\frac{3}{2}+\frac{1}{8}+\varepsilon}}{Q^{\frac{1}{2}}}
+\frac{X^{2+\frac{1}{8}+\varepsilon}}{Q^{\frac{3}{2}}}\nonumber\\
&\ll_{f,\varepsilon}&X^{\frac{3}{2}-\frac{1}{8}+\varepsilon}
\eea
uniformly for $1\leq h\leq X$.

Further, by (3.7), (3.13), Lemma 2 and Proposition 3.3, we have
\bna
\mathscr{D}_3(v,q,\beta)
&\ll_{f,\varepsilon}&
 \frac{1}{q}\frac{X}{1+|\beta|X}
\sum_{n_1|q}\sum_{\frac{n_1^2n_2}{q^3}X< X^{\varepsilon}(1+|\beta|X)^3}
\frac{\left|A_f(n_2,n_1)\right|}{n_1n_2}
\left|\Phi_{\beta}^+\left(\frac{n_1^2n_2}{q^3}\right)\right|\nonumber\\
&&\times\sum_{-\frac{3q}{2}<b\leq \frac{3q}{2}}|\Psi(b,q,\beta)|
|\mathscr{C}(0,0,b,n_1,n_2,h,v;q)|+1\nonumber\\
&\ll_{f,\varepsilon}&
\frac{X^{1+\varepsilon}}{1+|\beta|X}
\sum_{n_1\ll X^{\frac{1}{4}+\varepsilon}\atop n_1\equiv 0(\bmod q_1)}
\frac{q_1^2q_0^2q_3'^{\frac{3}{2}}(h,q_3')^{\frac{1}{2}}}{\sqrt{n_1}}\sum_{n_2=1}^{\infty}\frac{\left|A_f(n_2,n_1)\right|}{n_1n_2}
\left|\Phi_{\beta}^+\left(\frac{n_1^2n_2}{q^3}\right)\right|+1\nonumber\\
&\ll_{f,\varepsilon}&
X^{1+\varepsilon}(1+|\beta|X)
\sum_{n_1\ll X^{\frac{1}{4}+\varepsilon}\atop n_1\equiv 0(\bmod q_1)}\frac{q_1^2q_0^2q_3'^{\frac{3}{2}}
(h,q_3')^{\frac{1}{2}}}{\sqrt{n_1}}.
\ena
It follows from this estimate and (3.3) that
\bea
&&\sum_{q\leq Q}\int\limits_{|\beta|\leq \frac{1}{q Q}}\sum_{v\bmod q}
\varrho(v,q,\beta)\mathscr{D}_3(v,q,\beta)\mathrm{d}\beta\nonumber\\
&\ll_{f,\varepsilon}&X^{1+\varepsilon}\sum_{n_1\ll X^{\frac{1}{4}+\varepsilon}}n_1^{-\frac{1}{2}}
\sum_{q_1|n_1}\sum_{q_3'\leq Q/q_1 \atop q_3' \,\mathrm{square-free}}
\sum_{q_0\leq Q/(q_1q_3') \atop 4q_0 \, \mathrm{square-full}}
q_1^2q_0^2q_3'^{\frac{3}{2}}(h,q_3')^{\frac{1}{2}}
\int\limits_{|\beta|\leq \frac{1}{q_1q_0q_3' Q}}
(1+|\beta|X)\mathrm{d}\beta\nonumber\\
&\ll_{f,\varepsilon}&X^{1+\varepsilon}\sum_{n_1\ll X^{\frac{1}{4}+\varepsilon}}n_1^{-\frac{1}{2}}
\sum_{q_1|n_1}\sum_{q_3'\leq Q/q_1 \atop q_3' \,\mathrm{square-free}}
\sum_{q_0\leq Q/(q_1q_3') \atop 4q_0 \, \mathrm{square-full}}
q_1^2q_0^2q_3'^{\frac{3}{2}}(h,q_3')^{\frac{1}{2}}
 \left(\frac{1}{q_1q_0q_3'Q}+\frac{X}{(q_1q_0q_3'Q)^2} \right)\nonumber\\
&\ll_{f,\varepsilon}&\frac{X^{1+\varepsilon}}{Q}
\sum_{n_1\ll X^{\frac{1}{4}+\varepsilon}}n_1^{-\frac{1}{2}}
\sum_{q_1|n_1}q_1\sum_{q_3'\leq Q/q_1}(h,q_3')^{\frac{1}{2}}q_3'^{\frac{1}{2}}
\left(\frac{Q}{q_1q_3'}\right)^{\frac{3}{2}}\nonumber\\
&&+\frac{X^{2+\varepsilon}}{Q^2}
\sum_{n_1\ll X^{\frac{1}{4}+\varepsilon}}n_1^{-\frac{1}{2}}
\sum_{q_1|n_1}\sum_{q_3'\leq Q/q_1}(h,q_3')^{\frac{1}{2}}q_3'^{-\frac{1}{2}}
\left(\frac{Q}{q_1q_3'}\right)^{\frac{1}{2}}\nonumber\\
&\ll_{f,\varepsilon}&X^{1+\frac{1}{8}+\varepsilon}Q^{\frac{1}{2}}
+\frac{X^{2+\frac{1}{8}+\varepsilon}}{Q^{\frac{3}{2}}}\nonumber\\
&\ll_{f,\varepsilon}&X^{\frac{3}{2}-\frac{1}{8}+\varepsilon}
\eea
uniformly for $1\leq h\leq X$.

Moreover, by (3.7), (3.14), Lemma 2.2 and Propositions 3.3,
\bna
\mathscr{D}_5(v,q,\beta)
&\ll_{f,\varepsilon}&
\left(\frac{X}{1+|\beta|X}\right)^{\frac{1}{2}}
\sum_{n_1|q}\sum_{\frac{n_1^2n_2}{q^3}X< X^{\varepsilon}(1+|\beta|X)^3}\frac{\left|A_f(n_2,n_1)\right|}{n_1n_2}
\left|\Phi_{\beta}^+\left(\frac{n_1^2n_2}{q^3}\right)\right|\nonumber\\
&&\times\sum_{-\frac{3q}{2}<b_j\leq \frac{3q}{2}\atop j=1,2}|\Psi(b_1,q,\beta)|\Psi(b_2,q,\beta)|
|\mathscr{C}(0,b_1,b_2,n_1,n_2,h,v;q)|+1\nonumber\\
&\ll_{f,\varepsilon}&X^{\varepsilon}\left(\frac{X}{1+|\beta|X}\right)^{\frac{1}{2}}
\sum_{n_1\ll X^{\frac{1}{4}+\varepsilon}\atop n_1\equiv 0(\bmod q_1)}n_1^{-\frac{1}{2}}
q_1^3q_0^3q_3'^{\frac{5}{2}}(h,q_3')^{\frac{1}{2}}
\sum_{n_2=1}^{\infty}\frac{\left|A_f(n_2,n_1)\right|}{n_1n_2}
\left|\Phi_{\beta}^+\left(\frac{n_1^2n_2}{q^3}\right)\right|+1\nonumber\\
&\ll_{f,\varepsilon}&X^{\frac{1}{2}+\varepsilon}\left({1+|\beta|X}\right)^{\frac{3}{2}}
\sum_{n_1\ll X^{\frac{1}{4}+\varepsilon}\atop n_1\equiv 0(\bmod q_1)}n_1^{-\frac{1}{2}}
q_1^3q_0^3q_3'^{\frac{5}{2}}(h,q_3')^{\frac{1}{2}}.
\ena
Applying (3.3) again, we obtain
\bea
&&\sum_{q\leq Q}\int\limits_{|\beta|\leq \frac{1}{q Q}}\sum_{v\bmod q}
\varrho(v,q,\beta)\mathscr{D}_5(v,q,\beta)\mathrm{d}\beta\nonumber\\
&\ll_{f,\varepsilon}&X^{\frac{1}{2}+\varepsilon}
\sum_{n_1\ll X^{\frac{1}{4}+\varepsilon}}n_1^{-\frac{1}{2}}
\sum_{q_1|n_1}\sum_{q_3'\leq Q/q_1 \atop q_3' \,\mathrm{square-free}}
\sum_{q_0\leq Q/(q_1q_3') \atop 4q_0 \, \mathrm{square-full}}
q_1^3q_0^3q_3'^{\frac{5}{2}}(h,q_3')^{\frac{1}{2}}
\int\limits_{|\beta|\leq \frac{1}{q_1q_0q_3' Q}}
(1+|\beta|X)^{\frac{3}{2}}\mathrm{d}\beta\nonumber\\
&\ll_{f,\varepsilon}&X^{\frac{1}{2}+\varepsilon}\sum_{n_1\ll X^{\frac{1}{4}+\varepsilon}}n_1^{-\frac{1}{2}}
\sum_{q_1|n_1}\sum_{q_3'\leq Q/q_1 \atop q_3' \,\mathrm{square-free}}
\sum_{q_0\leq Q/(q_1q_3') \atop 4q_0 \, \mathrm{square-full}}
q_1^3q_0^3q_3'^{\frac{5}{2}}(h,q_3')^{\frac{1}{2}}
 \left(\frac{1}{q_1q_0q_3'Q}+\frac{X^{\frac{3}{2}}}{(q_1q_0q_3'Q)^{\frac{5}{2}}} \right)\nonumber\\
&\ll_{f,\varepsilon}&\frac{X^{\frac{1}{2}+\varepsilon}}{Q}
\sum_{n_1\ll X^{\frac{1}{4}+\varepsilon}}n_1^{-\frac{1}{2}}
\sum_{q_1|n_1}q_1^2\sum_{q_3'\leq Q/q_1}(h,q_3')^{\frac{1}{2}}q_3'^{\frac{3}{2}}
\left(\frac{Q}{q_1q_3'}\right)^{\frac{5}{2}}\nonumber\\
&&+\frac{X^{2+\varepsilon}}{Q^{\frac{5}{2}}}
\sum_{n_1\ll X^{\frac{1}{4}+\varepsilon}}n_1^{-\frac{1}{2}}
\sum_{q_1|n_1}q_1^{\frac{1}{2}}\sum_{q_3'\leq Q/q_1}(h,q_3')^{\frac{1}{2}}\frac{Q}{q_1q_3'}\nonumber\\
&\ll_{f,\varepsilon}&X^{\frac{1}{2}+\frac{1}{8}+\varepsilon}Q^{\frac{3}{2}}
+\frac{X^{2+\frac{1}{8}+\varepsilon}}{Q^{\frac{3}{2}}}\nonumber\\
&\ll_{f,\varepsilon}&X^{\frac{3}{2}-\frac{1}{8}+\varepsilon}
\eea
uniformly for $1\leq h\leq X$.

Lastly, by (3.7), (3.15), Lemma 2.2 and Propositions 3.3, we have
\bna
\mathscr{D}_7(v,q,\beta)
&\ll_{f,\varepsilon}& q
\sum_{n_1|q}\sum_{\frac{n_1^2n_2}{q^3}X< X^{\varepsilon}(1+|\beta|X)^3}\frac{\left|A_f(n_2,n_1)\right|}{n_1n_2}
\left|\Phi_{\beta}^+\left(\frac{n_1^2n_2}{q^3}\right)\right|\nonumber\\
&&\times\sum_{-\frac{3q}{2}<b_j\leq \frac{3q}{2}\atop 1\leq j\leq 3}|\Psi(b_1,q,\beta)||\Psi(b_2,q,\beta)|
|\Psi(b_3,q,\beta)|
|\mathscr{C}(b_1,b_2,b_3,n_1,n_2,h,v;q)|+1\nonumber\\
&\ll_{f,\varepsilon}&X^{\varepsilon}
\sum_{n_1\ll X^{\frac{1}{4}+\varepsilon}\atop n_1\equiv 0(\bmod q_1)}n_1^{-\frac{1}{2}}q_1^4q_0^4q_3'^{\frac{7}{2}}(h,q_3')^{\frac{1}{2}}\sum_{n_2=1}^{\infty}\frac{\left|A_f(n_2,n_1)\right|}{n_1n_2}
\left|\Phi_{\beta}^+\left(\frac{n_1^2n_2}{q^3}\right)\right|+1\nonumber\\
&\ll_{f,\varepsilon}&X^{\varepsilon}(1+|\beta|X)^2
\sum_{n_1\ll X^{\frac{1}{4}+\varepsilon}\atop n_1\equiv 0(\bmod q_1)}n_1^{-\frac{1}{2}}
q_1^4q_0^4q_3'^{\frac{7}{2}}(h,q_3')^{\frac{1}{2}}.
\ena
It follows that
\bea
&&\sum_{q\leq Q}\int\limits_{|\beta|\leq \frac{1}{q Q}}\sum_{v\bmod q}
\varrho(v,q,\beta)\mathscr{D}_7(v,q,\beta)\mathrm{d}\beta\nonumber\\
&\ll_{f,\varepsilon}&X^{\varepsilon}
\sum_{n_1\ll X^{\frac{1}{4}+\varepsilon}}n_1^{-\frac{1}{2}}
\sum_{q_1|n_1}\sum_{q_3'\leq Q/q_1 \atop q_3' \,\mathrm{square-free}}
\sum_{q_0\leq Q/(q_1q_3') \atop 4q_0 \, \mathrm{square-full}}
q_1^4q_0^4q_3'^{\frac{7}{2}}(h,q_3')^{\frac{1}{2}}
\int\limits_{|\beta|\leq \frac{1}{q_1q_0q_3' Q}}
(1+|\beta|X)^2\mathrm{d}\beta\nonumber\\
&\ll_{f,\varepsilon}&X^{\varepsilon}
\sum_{n_1\ll X^{\frac{1}{4}+\varepsilon}}n_1^{-\frac{1}{2}}
\sum_{q_1|n_1}\sum_{q_3'\leq Q/q_1 \atop q_3' \,\mathrm{square-free}}
\sum_{q_0\leq Q/(q_1q_3') \atop 4q_0 \, \mathrm{square-full}}
q_1^4q_0^4q_3'^{\frac{7}{2}}(h,q_3')^{\frac{1}{2}}
 \left(\frac{1}{q_1q_0q_3'Q}+\frac{X^2}{(q_1q_0q_3'Q)^3} \right)\nonumber\\
&\ll_{f,\varepsilon}&\frac{X^{\varepsilon}}{Q}
\sum_{n_1\ll X^{\frac{1}{4}+\varepsilon}}n_1^{-\frac{1}{2}}
\sum_{q_1|n_1}q_1^3\sum_{q_3'\leq Q/q_1}(h,q_3')^{\frac{1}{2}}q_3'^{\frac{5}{2}}
\left(\frac{Q}{q_1q_3'}\right)^{\frac{7}{2}}\nonumber\\
&&+\frac{X^{2+\varepsilon}}{Q^3}
\sum_{n_1\ll X^{\frac{1}{4}+\varepsilon}}n_1^{-\frac{1}{2}}
\sum_{q_1|n_1}q_1\sum_{q_3'\leq Q/q_1}(h,q_3')^{\frac{1}{2}}q_3'^{\frac{1}{2}}\left(\frac{Q}{q_1q_3'}\right)^{\frac{3}{2}}\nonumber\\
&\ll_{f,\varepsilon}&X^{\frac{1}{8}+\varepsilon}Q^{\frac{5}{2}}
+\frac{X^{2+\frac{1}{8}+\varepsilon}}{Q^{\frac{3}{2}}}\nonumber\\
&\ll_{f,\varepsilon}&X^{\frac{3}{2}-\frac{1}{8}+\varepsilon}
\eea
uniformly for $1\leq h\leq X$. By (3.2), (3.11) and (3.20)-(3.23), Theorem 1.1 follows.

\medskip

\section{Proof of Proposition 3.2}
\setcounter{equation}{0}
\medskip

We only prove (3.17) and (3.18) can be proved similarly.
Recall $\Phi_{\beta}^{\pm}(x)$ in (3.9) which we relabel as
\bea
\Phi_{\beta}^{\pm}(x)=\Phi_0(x,\beta)\pm\frac{1}{i\pi^3x}\Phi_1(x,\beta),
\eea
where for $\sigma>\max\limits_{1\leq j\leq 3}\{-1-\mathrm{Re}(\mu_j)-2k\}$,
\bea
\Phi_k(x,\beta)=\int\limits_{\mathrm{Re}(s)=\sigma}(\pi^3 x)^{-s}\prod_{j=1}^3
\frac{\Gamma\left(\frac{1+s+\mu_j+2k}{2}\right)}
{\Gamma\left(\frac{-s-\mu_j}{2}\right)}
\widetilde{\phi_{\beta}}(-s-k)\mathrm{d}s
\eea
with $\phi_{\beta}(x)=\phi\left(\frac{x-h}{X}\right)e(-\beta x)$. Note that
\bna
\phi_{\beta}^{(j)}(x)\ll_j \left(\frac{1+|\beta|X}{X}\right)^j.
\ena
By Lemma 2.2, we have
\bea
\sum_{n_2=1}^{\infty}\frac{|A_f(n_2,n_1)|}{n_1n_2}
\left|\Phi_{\beta}^{\pm}\left(\frac{n_1^2n_2}{q^3}\right)\right|=
\sum_{\frac{n_1^2n_2}{q^3}X< X^{\varepsilon}(1+|\beta|X)^3}\frac{|A_f(n_2,n_1)|}{n_1n_2}
\left|\Phi_{\beta}^{\pm}\left(\frac{n_1^2n_2}{q^3}\right)\right|+O_{f,\varepsilon}(1).
\eea
By (2.1) and (2.9), we have
\bea
\sum_{\frac{n_1^2n_2}{q^3}X\leq X^{\varepsilon}}\frac{|A_f(n_2,n_1)|}{n_2}
\left|\Phi_{\beta}^{\pm}\left(\frac{n_1^2n_2}{q^3}\right)\right|
&\ll_{f,\varepsilon}& X^{\varepsilon}(1+|\beta|X)
\max_{1\leq T\leq \frac{q^3X^{\varepsilon}}{n_1^2X}}\frac{1}{T}
\sum_{T\leq n_2\leq 2T}|A_f(n_2,n_1)|\nonumber\\
&\ll_{f,\varepsilon}& X^{\varepsilon}n_1(1+|\beta|X).
\eea

For $xX>X^{\varepsilon}$, by Lemma 2.3, we have
\bna
\Phi_k(x,\beta)&=&(\pi^3 x)^{k+1}\sum_{j=1}^\ell \int_0^{\infty}
\phi\left(\frac{u-h}{X}\right)e(-\beta u)
\left(a_k(j)e\left(3(x u)^{\frac{1}{3}}\right)+b_k(j)e\left(-3(x u)^{\frac{1}{3}}\right)\right)
\frac{\mathrm{d}u}{(\pi^3xu)^{\frac{j}{3}}}\\
&&+O_{f,\varepsilon,\ell}
\left((\pi^3x)^k(\pi^3xX)^{-\frac{\ell}{3}+\frac{1}{2}+\varepsilon}\right),
\ena
where $a_k(j)$, $b_k(j)$ are constants, and by (4.1),
\bea
\Phi_{\beta}^{\pm}(x)\ll_{f,\varepsilon,\ell}x\sum_{j=1}^\ell x^{-\frac{j}{3}}\left(|\mathcal{I}_j(x,\beta)|
+|\mathcal{J}_j(x,\beta)|\right)+(xX)^{-\frac{\ell}{3}+\frac{1}{2}+\varepsilon},
\eea
with
\bna
\mathcal{I}_j(x,\beta)&=&\int_0^{\infty}u^{-\frac{j}{3}}\phi\left(\frac{u-h}{X}\right)e(-\beta u)
e\left(3(xu)^{\frac{1}{3}}\right)\mathrm{d}u,\\
\mathcal{J}_j(x,\beta)&=&\int_0^{\infty}u^{-\frac{j}{3}}\phi\left(\frac{u-h}{X}\right)e(-\beta u)
e\left(-3(xu)^{\frac{1}{3}}\right)\mathrm{d}u.
\ena

By partial integration twice, we have
\bea
\mathcal{I}_j(x,\beta)&=&\frac{1}{\left(2\pi i x^{\frac{1}{3}}\right)^2}\int_0^{\infty}
\left(u^{\frac{2}{3}}\left(u^{\frac{2-j}{3}}\phi\left(\frac{u-h}{X}\right)e(-\beta u)\right)'\right)'
e\left(3(xu)^{\frac{1}{3}}\right)\mathrm{d}u\nonumber\\
&\ll & (xX)^{-\frac{2}{3}}X^{1-\frac{j}{3}}(1+|\beta|X)^2.
\eea
Similarly,
\bea
\mathcal{J}_j(x,\beta)
&=&\frac{1}{\left(2\pi i x^{\frac{1}{3}}\right)^2}\int_0^{\infty}
\left(u^{\frac{2}{3}}\left(u^{\frac{2-j}{3}}\phi\left(\frac{u-h}{X}\right)e(-\beta u)\right)'\right)'
e\left(-3(xu)^{\frac{1}{3}}\right)\mathrm{d}u\nonumber\\
&\ll & (xX)^{-\frac{2}{3}}X^{1-\frac{j}{3}}(1+|\beta|X)^2.
\eea
Taking $\ell=3$. By (4.5)-(4.7), we have
\bna
\Phi_{\beta}^{\pm}(x)\ll_{f,\varepsilon,\ell}(xX)^{\frac{1}{3}}(1+|\beta|X)^2
\sum_{j=1}^\ell (xX)^{-\frac{j}{3}}+(xX)^{-\frac{\ell}{3}+\frac{1}{2}+\varepsilon}
\ll_{f,\varepsilon} (1+|\beta|X)^2.
\ena
This estimate combined with (2.1) yields that
\bea
&&\sum_{X^{\varepsilon}<\frac{n_1^2n_2}{q^3}X< X^{\varepsilon}(1+|\beta|X)^3}\frac{|A_f(n_2,n_1)|}{n_2}
\left|\Phi_{\beta}^{\pm}\left(\frac{n_1^2n_2}{q^3}\right)\right|\nonumber\\
&\ll_{f,\varepsilon}&(1+|\beta|X)^2 (\log X)
\max_{\frac{q^3X^{\varepsilon}}{n_1^2X}<T< \frac{q^3X^{\varepsilon}(1+|\beta|X)^3}{n_1^2X}}\frac{1}{T}
\sum_{T\leq n_2\leq 2T}|A_f(n_2,n_1)|\nonumber\\
&\ll_{f,\varepsilon}&X^{\varepsilon}n_1(1+|\beta|X)^2.
\eea
Then Proposition 3.2 follows from (4.3), (4.4) and (4.8).

\medskip

\section{Estimation of the character sum $\mathscr{C}(b_1,b_2,b_3,n_1,n_2,h,v;q)$}
\setcounter{equation}{0}
\medskip

Let $b_1,b_2,b_3,n_1, n_2, v\in \mathbb{Z}$, $n_1|q$ and $1\leq h\leq X$.
We recall $\mathscr{C}(b_1,b_2,b_3,n_1,n_2,h,v;q)$ in (3.16) which we relabel as
\bea
\mathscr{C}(b_1,b_2,b_3,n_1,n_2,h,v;q)=\sideset{}{^*}\sum_{a\bmod q}e\left(\frac{ah-\overline{a}v}{q}\right)
G(a,b_1;q)G(a,b_2;q)G(a,b_3;q) S\left(-\overline{a},n_2;\frac{q}{n_1}\right),
\nonumber\\
\eea
where $S(m,n;c)$ is the classical Kloosterman sum and $G(a,b;q)$ is the Gauss sum in (3.5)
$$
G(a,b;q)=\sum_{x\bmod q}e\left(\frac{ax^2+bx}{q}\right).
$$

In this section, we shall prove Proposition 3.3. We need the following results for $G(a,b;q)$
(see Lemma 5.4.5 in \cite{Huxley})

\medskip

\noindent {\bf Lemma 5.1.} {\it (1) If $(q_1,q_2)=1$, then $G(a,b;q_1q_2)=G(aq_2,b;q_1)G(aq_1,b;q_2)$.

(2) For $(a,q)=1$, the sum $G(a,b;q)$ has absolute value $\sqrt{q}$ if $q$ is odd,
$\sqrt{2q}$ if $q$ is even, $q=2r$ and $ar+b$ is even.

(3) For $(2a,q)=1$, we have
\bna
G(a,0;q)=\left(\frac{a}{q}\right)\epsilon_q\sqrt{q},
\ena
where
$
\epsilon_q=\left\{\begin{array}{ll}
1,&\mbox{if $q\equiv 1(\bmod 4)$},\\
i,&\mbox{if $q\equiv -1(\bmod 4)$.}
\end{array}
\right.
$

(4) For $q$ odd, we have
\bna
G(a,b;q)=
e\left(-\frac{\bar{4}\bar{a}b^2}{q}\right)G(a,0;q).
\ena
}

In the estimate for $\mathscr{C}(b_1,b_2,b_3,n_1,n_2,h,v;q)$,
the main case of interest
is when $q$ is square-free and $n_1=1$. Thus we first extract from $q$ the
largest part related to $n_1$, the contribution of which to
$\mathscr{C}(b_1,b_2,b_3,n_1,n_2,h,v;q)$ will be trivially estimated
(using Weil's bound for Kloosterman sums). More precisely,
let $q_1$ be the largest factor of $q$ such that $q_1|n_1$
and $(q_1, q/q_1)=1$.
Let $q_2$ be the largest factor of $q/q_1$ such that $q_2|n_1^{\infty}$ and
$\left(q_2,\frac{q}{q_1q_2}\right)=1$.
Note that $q_2$ is square-full and $n_1|q_1q_2$. Denote
temporarily $q'=q_1q_2$ and $\widehat{q}=\frac{q'}{n_1}$.
Let $q_3=\frac{q}{q_1q_2}$.
Then any reduced residue class $a \bmod q$ can be written as
$a=a_1q_3+a_2q'$ with inverse $\overline{a}=\overline{a}_1q_3\overline{q}_3^2
+\overline{a}_2q'\overline{q'}^2$, where $a_1 \bmod q'$
and $a_2 \bmod q_3$ are some reduced residue classes modulo $q'$ and $q_3$, respectively.
Then by Lemma 5.1 (1), we have
\bea
&&\mathscr{C}(b_1,b_2,b_3,n_1,n_2,h,v;q)\nonumber\\
&=&\sideset{}{^*}\sum_{a_1\bmod q'}
e\left(\frac{a_1h-\overline{a}_1\overline{q}_3^2v}{q'}\right)
G(a_1q_3^2,b_1;q')G(a_1q_3^2,b_2;q')G(a_1q_3^2,b_3;q')
S(-\overline{a_1}\overline{q_3},n_2\overline{q_3}^2;\widehat{q})\nonumber\\
&&\times \sideset{}{^*}\sum_{a_2\bmod q_3}
e\left(\frac{a_2h-\overline{a}_2\overline{q'}^2v}{q_3}\right)
G(a_2q'^2,b_1;q_3)G(a_2q'^2,b_2;q_3)G(a_2q'^2,b_3;q_3)
S(-\overline{a_2}\overline{q'},n_2\overline{\widehat{q}}^2;q_3)
\nonumber\\
&=&\sideset{}{^*}\sum_{a_1\bmod q'}
e\left(\frac{a_1\overline{q}_3^2h-\overline{a_1}v}{q'}\right)
G(a_1,b_1;q')G(a_1,b_2;q')G(a_1,b_3;q')
S(-\overline{a_1}q_3,n_2\overline{q_3}^2;\widehat{q})\nonumber\\
&&\times \sideset{}{^*}\sum_{a_2\bmod q_3}
e\left(\frac{a_2\overline{q'}^2h-\overline{a}_2v}{q_3}\right)
G(a_2,b_1;q_3)G(a_2,b_2;q_3)G(a_2,b_3;q_3)
S(-\overline{a_2}q',n_2\overline{\widehat{q}}^2;q_3)
\nonumber\\
&:=&\mathscr{C}^*(b_1,b_2,b_3,n_1,n_2,h,v;q')
\mathscr{C}^{**}(b_1,b_2,b_3,n_1,n_2,h,v;q_3)
\eea
say.

By Lemma 5.1 (2), we have $G(a,b;q)\ll \sqrt{q}$. This estimate combined with Weil's bound
for Kloosterman sum gives
\bea
\mathscr{C}^*(b_1,b_2,b_3,n_1,n_2,h,v;q')\ll q'^{\frac{5}{2}}
\left(\overline{a_1}q_3,n_2\overline{q_3}^2,\frac{q'}{n_1}\right)^{\frac{1}{2}}
\left(\frac{q'}{n_1}\right)^{\frac{1}{2}}\tau\left(\frac{q'}{n_1}\right)
\ll \frac{q_1^3q_2^3\tau(q_1q_2)}{\sqrt{n_1}}.
\eea

Next we extract the
square-full part from the remaining part of $q$, that is to say,
we further factor
$q_3$ as $q_3=q_3'q_3''$ with $(q_3',2q_3'')=1$, $q_3'$ square-free and $4q_3''$ square-full.
Then
\bea
\mathscr{C}^{**}(b_1,b_2,b_3,n_1,n_2,h,v;q_3)
=\mathscr{C}^{**}_1(b_1,b_2,b_3,n_1,n_2,h,v;q'_3)
\mathscr{C}^{**}_2(b_1,b_2,b_3,n_1,n_2,h,v;q''_3),
\eea
where
\bna
\mathscr{C}^{**}_1(b_1,b_2,b_3,n_1,n_2,h,v;q'_3)&=&
\sideset{}{^*}\sum_{\gamma \bmod q_3'}
e\left(\frac{\gamma\overline{q'}^2\overline{q_3''}^2h-
\overline{\gamma}v}{q_3'}\right)
G(\gamma,b_1;q_3')G(\gamma,b_2;q_3')\\
&&
G(\gamma,b_3;q_3')
S(-\overline{\gamma}q'q_3'',n_2\overline{\widehat{q}}^2\overline{q_3''}^2;q_3')\\
\mathscr{C}^{**}_2(b_1,b_2,b_3,n_1,n_2,h,v;q''_3)
&=&\sideset{}{^*}\sum_{\gamma \bmod q_3''}
e\left(\frac{\gamma\overline{q'}^2\overline{q_3'}^2h-
\overline{\gamma}v}{q_3''}\right)
G(\gamma,b_1;q_3'')G(\gamma,b_2;q_3'')\\
&&
G(\gamma,b_3;q_3'')
S(-\overline{\gamma}q'q_3',n_2\overline{\widehat{q}}^2\overline{q_3'}^2;q_3'').
\ena
By Lemma 5.1 (2) and Weil's bound for Kloosterman sum, we have
\bea
\mathscr{C}^{**}_2(b_1,b_2,b_3,n_1,n_2,h,v;q''_3)
\ll q''^{3}\tau(q_3'').
\eea

To estimate $\mathscr{C}^{**}_1:=\mathscr{C}^{**}_1(b_1,b_2,b_3,n_1,n_2,h,v;q'_3)$, we factor $q_3'$ as $q_3'=p_1p_2\cdot\cdot\cdot p_s$, $p_i$ prime. Then
\bea
\mathscr{C}^{**}_1=\prod_{i=1}^s\mathscr{T}(b_1,b_2,b_3,\overline{q'}^2\overline{q_3''}^2\overline{p_i'}^2h,
q'q_3''p_i',n_2\overline{\widehat{q}}^2\overline{q_3''}^2\overline{p_i'}^2;p_i),
\eea
where $p_i'=q_3'/p_i$ and
\bna
\mathscr{T}(b_1,b_2,b_3,r_1h,r_2,r_3n_2;p)
=\sideset{}{^*}\sum_{x\bmod p}e\left(\frac{r_1hx-v\overline{x}}{p}\right)
G(x,b_1;p)G(x,b_2;p)G(x,b_3;p)S(-r_2\overline{x},r_3n_2;p)
\ena
with $(p,2)=1$ and $(r_i,p)=1$, i=1,2,3.

By Lemma 5.1 (3), (4), we have
\bna
G(x,b_j;p)=e\left(-\frac{\overline{4}\overline{x}b_j^2}{p}\right)G(x,0;p)
=e\left(-\frac{\overline{4}\overline{x}b_j^2}{p}\right)\left(\frac{x}{p}\right)
\epsilon_p\sqrt{p}.
\ena
Hence
\bea
&&\mathscr{T}(b_1,b_2,b_3,r_1h,r_2,r_3n_2;p)\nonumber\\
&=&\epsilon_p^3p^{\frac{3}{2}}\sideset{}{^*}\sum_{x\bmod p}\left(\frac{x}{p}\right)
e\left(\frac{r_1hx-\overline{4}(4v+b_1^2+b_2^2+b_3^2)\overline{x}}{p}\right)
S(-r_2\overline{x},r_3n_2;p).
\eea
Denote
\bna
\widetilde{\mathscr{T}}(p)=\sideset{}{^*}\sum_{x\bmod p}\left(\frac{x}{p}\right)
e\left(\frac{r_1hx-\overline{4}(4v+b_1^2+b_2^2+b_3^2)\overline{x}}{p}\right)
S(-r_2\overline{x},r_3n_2;p).
\ena

By (5.1)-(5.7), Proposition 3.3 follows from the following lemma.

\medskip

\noindent {\bf Lemma 5.2.} {\it We have
\bna
\widetilde{\mathscr{T}}(p)\ll (h,p)^{\frac{1}{2}}p.
\ena
}

\noindent {\bf Proof.} If $p|h$, then by Weil's bound for Kloosterman sum, we have
\bea
\widetilde{\mathscr{T}}(p)\ll p^{\frac{3}{2}}.
\eea

If $p\nmid h$, $p|n_2$, then $S(-r_2\overline{x},r_3n_2;p)=-1$ and
\bea
\widetilde{\mathscr{T}}(p)=-\sideset{}{^*}\sum_{x\bmod p}\left(\frac{x}{p}\right)
e\left(\frac{r_1hx-\overline{4}(4v+b_1^2+b_2^2+b_3^2)\overline{x}}{p}\right)
\ll p^{\frac{1}{2}}
\eea
by the bound for Sali\'{e} sum (see \cite{Iwan}, Corollary 4.10).

If $p\nmid h$, $p\nmid n_2$, we open the Kloosterman sum to obtain
\bna
\widetilde{\mathscr{T}}(p)=\sum_{x,y\in \mathbb{F}_p^{\times}}\left(\frac{x}{p}\right)
e\left(\frac{r_1hx-\overline{4}(4v+b_1^2+b_2^2+b_3^2)\overline{x}-r_2\overline{x}y
+r_3n_2\overline{y}}{p}\right).
\ena
The square-root cancellation for such twisted character sums was established in the more general case
in \cite{AS} (see also \cite{Fu}). To apply their result to the special character sum
$\widetilde{\mathscr{T}}(p)$, we consider the Newton polyhedron $\Delta(f)$ of
$f(x,y)=r_1hx-\overline{4}(4v+b_1^2+b_2^2+b_3^2)x^{-1}-r_2x^{-1}y+r_3n_2y^{-1}\in
\mathbb{F}_{p}^{\times}[x,y,(xy)^{-1}]$.
If $p\nmid 4v+b_1^2+b_2^2+b_3^2$,
then $\Delta(f)$ is the quadrilateral in $\mathbb{R}^2$ with vertices
$(1,0)$, $(-1,0)$, $(-1,1)$ and $(0,-1)$. If $p\mid 4v+b_1^2+b_2^2+b_3^2$,
then $\Delta(f)$ is the triangle in $\mathbb{R}^2$ with vertices
$(1,0)$, $(-1,1)$ and $(0,-1)$. Thus in both cases $\mathrm{dim}\Delta(f)=2$ and
$(0,0)$ is an interior point of $\Delta(f)$. Moreover, for $p\nmid 4v+b_1^2+b_2^2+b_3^2$, each of the following eight
polynomials
\bna
f_{\sigma}(x,y)&=&r_1hx, -\overline{4}(4v+b_1^2+b_2^2+b_3^2)x^{-1},
-r_2x^{-1}y,r_3n_2y^{-1}, r_1hx-r_2x^{-1}y, r_1hx+r_3n_2y^{-1},\\
&&-\overline{4}(4v+b_1^2+b_2^2+b_3^2)x^{-1}-r_2x^{-1}y,
-\overline{4}(4v+b_1^2+b_2^2+b_3^2)x^{-1}+r_3n_2y^{-1}
\ena
 and for $p\mid 4v+b_1^2+b_2^2+b_3^2$, each of the following six polynomials
\bna
f_{\sigma}(x,y)=r_1hx,
-r_2x^{-1}y,r_3n_2y^{-1}, r_1hx-r_2x^{-1}y, r_1hx+r_3n_2y^{-1},
-r_2x^{-1}y+r_3n_2y^{-1}
\ena
corresponding to the faces of $\Delta(f)$ not containing $(0,0)$,
the locus of
\bna
\frac{\partial f_{\sigma}}{\partial x}=\frac{\partial f_{\sigma}}{\partial y}=0
\ena
is empty in $\left(\overline{F}_p^{\times}\right)^2$. Thus in both cases $f$ is nondegenerate
with respect to $\Delta(f)$.
By \cite{AS} or \cite{Fu}, we have
\bea
\widetilde{\mathscr{T}}(p)\ll p.
\eea
Then Lemma 5.2 follows from (5.8)-(5.10). \hfill $\Box$

\bigskip

\bigskip

\noindent
{\sc Acknowledgements.}
The author would like to express heartfelt thanks to the referees
for their important and useful comments.
This work is supported by
the National Natural Science Foundation of China (Grant No. 11101239) and
Young Scholars Program of Shandong University, Weihai (Grant No. 2015WHWLJH04).

\end{document}